\documentclass[12pt]{article}

\usepackage[margin=20truemm]{geometry}

\usepackage[
    uniquelist=false,
    style=alphabetic,
    maxnames=8,
    minnames=5,
    abbreviate=true,
    sorting=nty
]{biblatex}

\usepackage{latexsym}
\usepackage{amsmath}
\usepackage{amssymb}
\usepackage{amsthm}

\usepackage[all]{xy}

\usepackage{color}

\usepackage{amsfonts}
\usepackage[usenames,dvipsnames]{xcolor}

\usepackage{varioref}
\usepackage{hyperref}
\usepackage[noabbrev,capitalise,nameinlink]{cleveref}
\hypersetup{colorlinks={true},linkcolor={blue},citecolor=blue}

\def\C{\mathbb{C}}
\def\Q{\mathbb{Q}}
\def\Z{\mathbb{Z}}
\def\F{\mathbb{F}}
\def\T{\mathbb{T}}
\def\p{\mathfrak{p}}

\def\O{\mathcal{O}}
\def\SL{\text{SL}} 
\def\sl{\text{sl}}
\def\tr{\text{tr}}
\def\Ad{{\rm Ad}}
\def\Int{\text{int}} 
\def\phi{\varphi} 
\def\Coker{\text{Coker}}
\def\del{\partial}
\def\Sel{{\rm Sel}}
\def\:{\colon}
\def\boldrho{{\boldsymbol \rho}}
\def\rhobar{\overline{\rho}}

\newcommand\blfootnote[1]{
  \begingroup
  \renewcommand\thefootnote{}\footnote{#1}
  \addtocounter{footnote}{-1}
  \endgroup
}

\begin{document}

\begin{center}
 {\Large \bf 
 On adjoint homological Selmer modules for ${\rm SL}_2$-representations 
 of knot groups}
 \end{center}

\begin{center}
Takahiro KITAYAMA, Masanori MORISHITA, Ryoto TANGE, and Yuji TERASHIMA
\end{center}



\noindent
{\bf Abstract.}
We introduce the adjoint homological Selmer module for 
an ${\rm SL}_2$-representation of a knot group, 
which may be seen as a
knot
theoretic 
analogue of the dual 
adjoint Selmer module for a Galois representation. 
We then show finitely generated torsion-ness 
of our adjoint Selmer module, 
which are widely known as conjectures in number theory, 
and give some concrete examples. 

\begin{center}
{\bf Introduction} 
\end{center}

In this paper, we continue our study on the interplay between knot theory and number theory, and 
introduce 
an adjoint Selmer module for the adjoint representation of a knot group representation. 
To pursue the analogy more strictly, 
we especially consider
the Fitting ideal of the adjoint Selmer module for 
the universal deformation of a knot group representation, 
which may be seen as an analogue of the algebraic $p$-adic $L$-function associated with the adjoint Selmer module for 
the universal deformation of a Galois representation. 
We then show finitely generated torsion-ness 
of our adjoint Selmer module, 
which are widely known as conjectures in number theory. 
Main ingredients of our proof are 
Euler characteristics and 
$\gamma$-regular 
representations. 

Let us recall the notion of a $\gamma$-regular representation in 3-dimensional
hyperbolic geometry (
\cite[3.3]{porti1997torsion}). 
For simplicity, we consider the case for a hyperbolic knot. 
Let $K$ be a hyperbolic knot in the 3-sphere $S^3$ and let
$E_K := S^3 \setminus \Int(V_K)$ be the knot exterior, where 
$V_K$ is a tubular neighborhood of $K$. 
Let $G_K := \pi_1(E_K)$ be the knot group and let 
$D_K := \pi_1(\partial E_K)$. 
Let $\rho \: G_K \to {\rm SL}_2(\C)$ be an irreducible representation such that 
$\rho|_{D_K}$ is hyperbolic or parabolic. 
Let ${\rm Ad}(\rho)$ be the adjoint representation of $\rho$, namely, 
the representation space is ${\rm sl}_2(\C)$ on which $G_K$ acts by 
${\rm Ad}(\rho)(X) := \rho(g) X \rho(g)^{-1}$ 
for $g \in G_K$ and $X \in {\rm sl}_2(\C)$. 
For a simple closed curve $\gamma$ on $\partial E_K$, 
let $I_{\gamma}$ be the subgroup of $D_K$ generated by $\gamma$, 
and we assume that $\rho(\gamma) \neq \pm I$. 
The representation $\rho$ is said to be $\gamma$-regular if the natural homomorphism

$$
\varphi_{\gamma} \: H_1(I_{\gamma}; {\rm Ad}(\rho)) \to H_1(G_K; {\rm Ad}(\rho) )
$$
is surjective. 

Let $X(G_K)$ be the character variety of $\SL_2(\C)$-representations of $G_K$ and
$X_0(G_K)$ be the irreducible component of $X(G_K)$ containing $\chi_{0} := {\rm tr}( \rho_0 )$, 
where $\rho_0$ is a lift of the holonomy representation of the complete hyperbolic 
structure on $\Int(E_K)$. 
Suppose $\chi_{\rho} := \tr(\rho) \in X_0(G_K)$. 
Joan Porti \cite{porti1997torsion} proved 
that when $\rho$ is $\gamma$-regular, 
the evaluation map
$$
X_0(G_K) \to \C; \ \chi \mapsto \chi(\gamma)
$$
gives an analytic isomorphism on a neighborhood of $\chi_\rho$. 
In particular, $\rho_0$ is $\gamma$-regular for any $\gamma (\neq 1) \in D_K$  
and the above parametrization around $\chi_0$ is
due to Thurston (\cite{thurston1979geometry}). 

Let $\overline{\Q}$ be the field of all algebraic numbers, 
$S$ be a finite set of primes, and 
$R$ be the ring of $S$-integers in some finite algebraic number field. 
Since $X_0(G_K)$ is defined over a finite algebraic number field (\cite{long2010fields}), 
$\overline{\Q}$-rational points are dense in $X_0(G_K)$ and such a point corresponds to a representation 
$\rho \: G_K \to \SL_2(R)$. 
For example, hyperbolic Dehn surgery points are such algebraic points. 
So it may be natural to ask the {\it integral} $\gamma$-regularity 
for $\rho \colon G_K \to \SL_2(R)$, 
namely, whether the natural $R$-homomorphism
$$
\phi_{\gamma, R} \: H_1(I_{\gamma}; \Ad(\rho)) \to H_1(G_K; \Ad(\rho) )
$$
is surjective or not. 
It amounts to study the $R$-module $\Coker(\phi_{\gamma, R})$. 
It turns out, from the viewpoint of arithmetic topology (\cite{morishita2012knots}), that 
the $R$-module $\Coker(\phi_{\gamma, R})$ may be regarded as an analogue of 
the dual adjoint Selmer module in Iwasawa theory for 
$p$-adic Galois representations 
(\cite{greenberg1989iwasawa}, \cite{hida2000adjoint}, \cite[Chapter 5]{hida2000modular}, \cite{hida2006hilbert}). 
See \cite{kitayama2018certain} for the analogy between twisted knot modules and dual Selmer modules.

We also give some concrete examples and 
guess a relation between the well-known topological invariants 
(Porti's torsions). 
Further relations between 
the Iwasawa class formula and Iwasawa invariants 
are expected, 
and generalizations 
such as 
symmetric power representations and quantum representations 
are also expected. 

This 
paper 
is organized as follows. 
In Section 1, 
we introduce the notion of an adjoint homological $\gamma$-Selmer module for 
any $\SL_2$-representation of $G_K$ over any integral domain $A$. 
In Section 2, 
we 
study its $A$-module structure for a Noetherian 
unique 
factorization domain $A$ 
and give a non-trivial example 
for the figure-eight knot. 
In Section 3, 
we give an example for 
a lift $\rho_0$ 
of the holonomy representation of 
the figure eight knot where $\rho_0$ is not 
$\mu$-regular 
over the number ring $\Z[(1 + \sqrt{-3})/2]$ 
for a meridian $\mu$. 
In Section 4, 
we consider the case of universal deformations, 
which are studied in 
\cite{kitayama2018certain}, 
\cite{morishita2017universal}, and  
\cite{tange2021non}. 
Lastly, 
in Section 5, 
by using 
the parameters of 
the universal deformation and 
the character variety of 
the one-dimensional representation of knot groups, 
namely, the parameter $t$ of the Alexander polynomial, 
we generalize the adjoint homological $\gamma$-Selmer module 
to the case of two-variable. 

\ 

\noindent
{\it Acknowledgements.}
We would like to thank 
Shinya Harada, 
Haruzo Hida, 
Alan Reid, 
Jun Ueki, 
Yoshikazu Yamaguchi, 
and the anonymous referee 
for helpful communications. 
We also thank Megumi Takata for his help with the proof of Proposition 4.1.
This work is
partially supported by 
JSPS KAKENHI Grant Numbers JP18K13404, JP18KK0071, JP18KK0380, JP21H00986, 
JP17H02837, 
JP17K05243, JP21K03240 
and by 
JST CREST Grant Number JPMJCR14D6.

\begin{center}
{\bf 1. Adjoint homological $\gamma$-Selmer modules} 
\end{center}

Let $K$ be a knot in the 3-sphere $S^3$. 
Let $V_K$ be a tubular neighborhood of $K$ and 
let $E_K := S^3 \setminus \Int(V_K)$ be the knot exterior. 
Let $G_K := \pi_1(E_K)$ be the knot group of $K$ and 
let $D_K := \pi_1(\partial V_K) = \pi_1(\partial E_K)$. 
For $\gamma \in D_K$, let $I_{\gamma}$ be the subgroup of $D_K$ generated by $\gamma$. 

Let $A$ be an integral domain with quotient field $Q(A)$. 
Let $\rho : G_K \rightarrow {\rm SL}_2(A)$ be a representation of $G_K$. 
Let $\Ad(\rho)$ denote the adjoint representation of $\rho$, namely, 
the representation space is the $A$-module $V_{\rho} := \sl_2(A)$ on which $G_K$ 
acts by 
$
\Ad(\rho(g)) (X) 
:= \rho(g) X \rho(g)^{-1}$ for $g \in G_K$ and $X \in {\rm sl}_2(A)$. 
For $\gamma \in D_K$, 
the natural homomorphism $I_{\gamma} \subset D_K \to G_K$ induces 
the $A$-homomorphism of group homology groups with coefficients in ${\rm Ad}(\rho)$: 
$$
\varphi_{\gamma} \: 
H_1(I_{\gamma}, \Ad(\rho)) 
\longrightarrow H_1(G_K, \Ad(\rho) ). 
$$
We then define the {\em adjoint homological $\gamma$-Selmer module}
\hspace*{-0.3cm}
\footnote{
Selmer modules in number theory are defined 
by using cohomology groups 
and 
so our homological Selmer modules 
are analogues of the dual of adjoint Selmer modules 
in number theory 
(cf. \cite{greenberg1989iwasawa}, \cite{hida2000modular}, \cite{hida2006hilbert}).
} 
attached to the representation $\rho$ by the $A$-module 
$$ {\rm Sel}^{\rm h}_{\gamma}({\rm Ad}(\rho)) := {\rm Coker}(\varphi_\gamma).$$

\noindent 
We say that 
$\rho$ is {\em $\gamma$-regular over $A$} 
if $\Sel^{h}_{\gamma}(\Ad(\rho)) = \{ 0 \}$. 
It is easy to see that 
${\rm Sel}^{\rm h}_{\gamma_1}({\rm Ad}(\rho))$ and 
${\rm Sel}^{\rm h}_{\gamma_2}({\rm Ad}(\rho))$ are 
$A$-isomorphic if 
$\gamma_1$ and $\gamma_2$ are conjugate in $G_K$. 
In particular, when $\gamma$ is a meridian 
of a knot $K$, 
we simply write 
${\rm Sel}^{\rm h}({\rm Ad}(\rho))$ 
for ${\rm Sel}^{\rm h}_{\gamma}({\rm Ad}(\rho))$ 
and call it the 
{\it adjoint homological Selmer module} attached to $\rho$.

Let $\rho_{Q(A)} : G_K \rightarrow {\rm SL}_2(Q(A))$ be the representation obtained from $\rho$ by the extension of scalars. \\
\\
{\bf Lemma 1.1.} {\em Suppose that ${\rm Sel}^{\rm h}_{\gamma}({\rm Ad}(\rho_{Q(A)})) = 0$. Then ${\rm Sel}^{\rm h}_{\gamma}({\rm Ad}(\rho))$ is a finitely generated torsion  $A$-module.}\\
\\
{\em Proof.} 
Since $G_K$ is a finitely presented group and ${\rm Ad}(\rho)$ is a free $A$-module of finite rank, $H_1(G_K, {\rm Ad}(\rho))$ is a finitely generated $A$-module and hence ${\rm Sel}^{\rm h}_{\gamma}({\rm Ad}(\rho))$ is so. So it suffices to show that ${\rm Sel}^{\rm h}_{\gamma}({\rm Ad}(\rho))$ is a torsion $A$-module. Let ${\rm tor}({\rm Sel}^{\rm h}_{\gamma}({\rm Ad}(\rho)))$ denote the torsion $A$-submodule of ${\rm Sel}^{\rm h}_{\gamma}({\rm Ad}(\rho))$. We note
$$ {\rm tor}({\rm Sel}^{\rm h}_{\gamma}({\rm Ad}(\rho))) = {\rm Ker}({\rm Sel}^{\rm h}_{\gamma}({\rm Ad}(\rho)) \longrightarrow {\rm Sel}^{\rm h}_{\gamma}({\rm Ad}(\rho))\otimes_A Q(A)). \leqno{(1.1.1)}$$

\noindent
By the definition of ${\rm Sel}^{\rm h}_{\gamma}({\rm Ad}(\rho))$, we have the exact sequence of $A$-modules
$$  H_1(I_{\gamma}, {\rm Ad}(\rho)) \stackrel{\varphi_{\gamma}}{\longrightarrow} H_1(G_K, {\rm Ad}(\rho)) \longrightarrow {\rm Sel}^{\rm h}_{\gamma}({\rm Ad}(\rho)) \longrightarrow 0. \leqno{(1.1.2)}$$

\noindent
Since $Q(A)$ is a flat $A$-module and 
$H_1(X, {\rm Ad}(\rho_{Q(A)})) = H_1(X, {\rm Ad}(\rho))\otimes_A Q(A)$ for 
$X = I_{\gamma}$, $G_K$, 
tensoring $Q(A)$ with (1.1.2) over $A$ 
and the assumption yield
$$ {\rm Sel}^{\rm h}_{\gamma}({\rm Ad}(\rho)) \otimes_A Q(A) = 
{\rm Sel}^{\rm h}_{\gamma}({\rm Ad}(\rho_{Q(A)})) = 0. \leqno{(1.1.3)}$$

\noindent
By (1.1.1) and (1.1.3), we have ${\rm tor}({\rm Sel}^{\rm h}_{\gamma}({\rm Ad}(\rho)))  = {\rm Sel}^{\rm h}_{\gamma}({\rm Ad}(\rho))$, and hence the assertion follows. 

$\hfill \Box$ 

\begin{center}
{\bf 2. 
Presentations 
of 
adjoint homological Selmer modules
} \\ 
\end{center}

We keep the same notations as in Section 1. 
In this Section, 
we assume that $A$ is a Noetherian 
unique 
factorization domain 
or 
a 
formal 
power series over a field. 
We 
study 
the $A$-module structure
of the adjoint homological Selmer module
${\rm Sel}^{\rm h}_{\gamma}({\rm Ad}(\rho))$. 

We take a Wirtinger presentation of $G_K$: 
$$G_K 
= \langle g_1,\dots , g_n \ | \ r_1 = \cdots = r_{n-1} = 1 \rangle, 
\leqno{(2.1)}$$

\noindent
where $g_1,\dots , g_n$ $(n \geq 2)$ represent 
meridians 
of a knot $K$, 
and let $\rho \: G_K \to \SL_2(A)$ be a representation. 
Let $F$ be the free group on the words $g_1, \dots , g_n$ and 
let $\pi \: A[F] \rightarrow A[G_K]$ be the natural homomorphism
of group rings. 
We write the same $g_i$ for the image of $g_i$ in $G_K$. 
Let $V_{\rho}$ be the representation space 
$$V_{\rho} := 
{\rm sl}_2(A) 
= A v_1 
\oplus A v_2 
\oplus A v_3, $$ 
where 
$$v_1 := \begin{pmatrix} 0 & 1 \\ 0 & 0 \end{pmatrix}, \; 
v_2 := \begin{pmatrix} 1 & 0 \\ 0 & -1 \end{pmatrix}, \; 
v_3 := \begin{pmatrix} 0 & 0 \\ 1 & 0 \end{pmatrix}. $$
In the following, 
taking $\gamma = g_1$, 
we shall compute 
a presentation matrix 
of 
${\rm Sel}^{\rm h}_{}({\rm Ad}(\rho)) = 
{\rm Sel}^{\rm h}_{\gamma}({\rm Ad}(\rho))$ 
over $A$, 
under some (mild) assumptions. \\

First, we compute 
$H_1( I_{\gamma}, \Ad(\rho) )$. 
Let $W_1$ be the CW complex attached to the presentation 
$\langle g_1 | - \rangle$, 
which is homotopically equivalent to the circle $S^1$. 
We consider the chain complex 
$
C_*(W_1, V_{\rho})$: 
$$
0 \longrightarrow 
C_1(W_1, V_{\rho}) \overset{d_1}{\longrightarrow} 
C_0(W_1, V_{\rho}) \longrightarrow 0, 
$$
defined by 
$C_1(W_1, V_{\rho}) := V_{\rho}$, 
$C_0(W_1, V_{\rho}) := V_{\rho}$ and 
$d_1 := {\rm Ad}({\rho})(g_1) - I$. 
So 
we have 
$$H_1( I_{\gamma}, {\rm Ad}({\rho}) ) 
= H_1(W_1,V_{\rho}) 
= {\rm Ker}(d_1) 
= \{ v \in V_{\rho} \ | \ {\rho}(g_1)v = v  {\rho}(g_1)\}. 
\leqno{(2.2)}$$ 

Next, we compute 
$H_1( G_K, {\rm Ad}({\rho}) )$. 
Let $W_2$ be the CW complex attached to 
the presentation (2.1). 
We note that the knot exterior $E_K$ 
and the CW complex $W_2$ are homotopically equivalent by Whitehead's theorem
because they are both the Eilenberg-MacLane space $K(G_K,1)$. 
We consider the chain complex 
$
C_*(W_2, V_{\rho})$: \\ 
$$
0 \longrightarrow 
C_2(W_2, V_{\rho}) \overset{\partial_2}{\longrightarrow} 
C_1(W_2, V_{\rho}) \overset{\partial_1}{\longrightarrow} C_0(W_2, V_{\rho}) \longrightarrow 0,
$$
defined by 
$$ \left\{ \begin{array}{l} 
C_2(W_2, V_{\rho}) := (V_{\rho})^{\oplus (n-1)},\\
C_1(W_2, V_{\rho}) := (V_{\rho})^{\oplus n},\\
C_0(W_2, V_{\rho}) := V_{\rho},
\end{array} \right.
\;\;
\left\{ \begin{array}{l}
\partial_2 := 
\begin{pmatrix}
{\rm Ad}({\rho}) \circ \pi \left( \frac{\partial r_i}{\partial g_j} \right)
\end{pmatrix}, \\ 
\partial_1 := 
\begin{pmatrix}
{\rm Ad}({\rho})(g_1) - I \\ 
\vdots \\ 
{\rm Ad}({\rho})(g_n) - I 
\end{pmatrix}, 
\end{array} \right.
$$
where $\frac{\partial}{\partial g_i} \: 
A[F] \rightarrow A[F]$ 
denotes the Fox derivative over $A$, 
extended from $\Z$ (\cite{fox1953free}), 
and $\partial_2$ is regarded as a (big) 
$(n - 1)\times n$ 
matrix whose 
$(i,j)$-entry is the $3 \times 3$ matrix 
${\rm Ad}({\rho}) 
\circ \pi \left( \frac{\partial r_i}{\partial g_j} \right)$. 
So we have 
$$
H_1( G_K, {\rm Ad}({\rho}) ) 
= H_1(W_2,V_{\rho}) 
= {\rm Ker}(\partial_1)/{\rm Im}(\partial_2). 
\leqno{(2.3)}$$ 

By (2.2), (2.3) and 
the definition of 
$\varphi_{\gamma}: H_1(I_{\gamma}, {\rm Ad}({\rho})) \rightarrow 
H_1(G_K, {\rm Ad}({\rho}) )$, 
we can regard 
$\varphi_{\gamma}$ 
as 
the 
composition of 
an 
inclusion map 
$$
{\rm Ker}(d_1) \hookrightarrow {\rm Ker}(\partial_1);  \; 
v \mapsto 
\begin{pmatrix} v \\ 0 \\ \vdots \\ 0 \end{pmatrix} 
$$
and the 
natural $A$-homomorphism 
${\rm Ker}(\partial_1) \rightarrow {\rm Ker}(\partial_1)/{\rm Im}(\partial_2)$. 
So we have the isomorphism of $A$-modules 
$${\rm Sel}^{\rm h}({\rm Ad}({\rho})) = {\rm Ker}(\partial_1)/({\rm Im}(\partial_2)+{\rm Ker}(d_1)), $$ 
and 
the exact sequence of $A$-modules: 
$$ 0 \rightarrow 
{\rm Sel}^{\rm h}({\rm Ad}({\rho})) \rightarrow 
V_{\rho}^{\oplus n} / ({\rm Im}(\partial_2)+{\rm Ker}(d_1)) \rightarrow 
V_{\rho}^{\oplus n} / {\rm Ker}(\partial_1) \rightarrow 
0. 
\leqno{(2.4)}$$

In the following, we make the assumptions 
$$ 
\left\{ \begin{array}{l} 
{\rm (A1)} \  {\rm Ker}(d_1) 
= A v_0 \ 
\textrm{for some} \ v_0 \in {\rm Ker}(d_1), \\
{\rm (A2)} \  \del_{1} \ \textrm{induces the} \ A \textrm{-isomorphism} \ 
V_{\rho}^{\oplus n} / {\rm Ker}(\partial_1) 
\simeq V_{\rho}, 
\end{array} \right. 
$$

\noindent 
Let 
$f : V_{\rho} \rightarrow V_{\rho}$ 
be the $A$-homomorphism 
defined by 
$$f(v_1) = v_0, \ f(v_2) = O, \ f(v_3) = O \ \ (O{\rm : zero \ matrix}). $$
Then we have 
$${\rm Ker}(d_1) = {\rm Im}(f). $$ 
Let 
$\delta : V_{\rho}^{\oplus n} = 
V_{\rho}^{\oplus (n-1) } \oplus V_{\rho} \rightarrow 
V_{\rho}^{\oplus n}$ 
be the $A$-homomorphism 
defined by 
$$\delta \begin{pmatrix} v \\ w \end{pmatrix} 
:= \partial_2(v) 
+ \begin{pmatrix} f(w) \\ 0 \\ \vdots \\ 0 \end{pmatrix}. $$ 
Then we have 
$${\rm Im}(\delta) = {\rm Im}(\partial_2) + {\rm Ker}(d_1). $$
By (2.4) and (A2), 
we have 
the exact sequence of $A$-modules  
$$0 \rightarrow 
{\rm Sel}^{\rm h}({\rm Ad}({\rho})) \rightarrow 
V_{\rho}^{\oplus n} / {\rm Im}(\delta) \rightarrow 
V_{\rho} \rightarrow 0. 
\leqno{(2.5)}$$

\noindent
We consider 
the $3n$-basis of $V_{\rho}^{\oplus n}$ \\ 
$$
{\bf v}_1^1 := \begin{pmatrix} v_1 \\ 0 \\ \vdots \\ 0 \end{pmatrix}, \ 
{\bf v}_2^1 := \begin{pmatrix} v_2 \\ 0 \\ \vdots \\ 0 \end{pmatrix}, \ 
{\bf v}_3^1 := \begin{pmatrix} v_3 \\ 0 \\ \vdots \\ 0 \end{pmatrix}, \  
$$
$$
{\bf v}_1^2 := \begin{pmatrix} 0 \\ v_1 \\ 0 \\ \vdots \\ 0 \end{pmatrix}, \ 
{\bf v}_2^2 := \begin{pmatrix} 0 \\ v_2 \\ 0 \\ \vdots \\ 0 \end{pmatrix}, \ 
{\bf v}_3^2 := \begin{pmatrix} 0 \\ v_3 \\ 0 \\ \vdots \\ 0 \end{pmatrix}, 
$$
$$\vdots$$ 
$$
{\bf v}_1^n := \begin{pmatrix} 0 \\ \vdots \\ 0 \\ v_1 \end{pmatrix}, \ 
{\bf v}_2^n := \begin{pmatrix} 0 \\ \vdots \\ 0 \\ v_2 \end{pmatrix}, \ 
{\bf v}_3^n := \begin{pmatrix} 0 \\ \vdots \\ 0 \\ v_3 \end{pmatrix}. 
$$ 
Then the presentation matrix $D$ of $\delta$ with respect to 
the basis 
$\{ {\bf v}_h^j \}$ 
($1 \leq h \leq 3$, $1 \leq j \leq n$) 
is given by the following form: 
$$D = \begin{pmatrix} 
A_{1,1} & A_{2,1} & \cdots & A_{n,1} \\
A_{1,2} & A_{2,2} & \cdots & A_{n,2} \\
\vdots & \vdots & \vdots & \vdots \\
A_{1,n-1} & A_{2,n-1} & \cdots & A_{n,n-1} \\ 
B & O & \cdots & O
\end{pmatrix}. 
\leqno{(2.6)}
$$

\noindent
Here, $A_{j,i}$ is the presentation matrix of 
${\rm Ad}({\rho}(\partial r_i/\partial g_j)) \: V_{\rho} \rightarrow V_{\rho}$ 
with respect to the basis $\{ v_1, v_2, v_3 \}$ 
and 
$B$ is the presentation matrix of $f : V_{\rho} \rightarrow V_{\rho}$ 
with respect to the basis $\{ v_1, v_2, v_3 \}$.

Summing up the above discussion, 
we have 
the following theorem. \\ 

\noindent 
{\bf Theorem 2.7.} {\em 
Let the notations and assumptions be as above. 
We have the exact sequence of $A$-modules 
$$ 
0 \rightarrow 
{\rm Sel}^{\rm h}({\rm Ad}({\rho})) \rightarrow 
V_{\rho}^{\oplus n} / {\rm Im}(\delta) \overset{\partial_1}{\longrightarrow} 
V_{\rho} \rightarrow 
0, 
$$
where a presentation matrix of 
$\delta \: V_{\rho}^{\oplus n} \rightarrow 
V_{\rho}^{\oplus n}$ 
is given by the matrix $D$ in (2.6), and 
for an integer $d \geq 0$, 
the $d$-th Fitting ideal of ${\rm Sel}^h(\Ad(\rho))$ over $A$ is 
generated by the greatest common divisor of all $(d + 3)$-minors of $D$. 
} \\ 

\noindent 
{\em Proof.} 
The latter assertion follows from that $V_{\rho}$ is a free $A$-module of rank $3$ 
(cf. \cite[7.2]{kawauchi2012survey}). 

\hfill $\Box$

\ 

\noindent 
{\bf Example 2.8.}
Let $K$ be the figure-eight knot, 
whose knot group is given by 
$$G_{K} = \langle g_1, g_2 \ | \ g_1 g_2^{-1} g_1^{-1} g_2 g_1 = g_2 g_1 g_2^{-1} g_1^{-1} g_2 \rangle. $$
We consider the case of $\gamma$ being the meridian 
$\mu = g_1$. 
Let 
$\rho_A \: G_{K} \to {\rm SL}_{2} ( A )$ 
be the representation 
given by the following 
(cf. \cite[Section 4]{kitayama2018certain}, \cite{riley1984nonabelian}): 
$$\rho_A(g_1) 
= \begin{pmatrix} 
\frac{x + \sqrt{x^{2} - 4}}{2} & 
1 \\ 
0 & 
\frac{x - \sqrt{x^{2} - 4}}{2} 
\end{pmatrix}, $$ 
$$\rho_A(g_2) 
= \begin{pmatrix} 
\frac{x + \sqrt{x^{2} - 4}}{2} & 
0 \\ 
- \{ x^{2} - y( x ) - 2 \} & 
\frac{x - \sqrt{x^{2} - 4}}{2} 
\end{pmatrix}, $$ 
where 
$A = \C[[ \sqrt{x-1} ]]$, and 
$y(x) =  \frac{x^2 + 1 + \sqrt{(x^2 - 1) (x^2 - 5)}}{2}$. 
We can rewrite these matrices as 
$$\rho_A(g_1) 
= \begin{pmatrix} 
\frac{(s^2 + 1) + \sqrt{(s^2 + 1)^{2} - 4}}{2} & 
1 \\ 
0 & 
\frac{(s^2 + 1) - \sqrt{(s^2 + 1)^{2} - 4}}{2} 
\end{pmatrix}, $$ 
$$\rho_A(g_2) 
= \begin{pmatrix} 
\frac{(s^2 + 1) + \sqrt{(s^2 + 1)^{2} - 4}}{2} & 
0 \\ 
- \{ (s^2 + 1)^{2} - y(s^2 + 1) - 2 \} & 
\frac{(s^2 + 1) - \sqrt{(s^2 + 1)^{2} - 4}}{2} 
\end{pmatrix}, $$ 
where we use 
$s := \sqrt{x - 1}$. 
Since 
$${\rm Ker}(\partial_1) 
= A 
\begin{pmatrix} 0 \\ * \\ 0 \\ 0 \\ 0 \\ 0 \end{pmatrix} \oplus  
A 
\begin{pmatrix} 0 \\  0 \\ * \\ 0 \\ 0 \\ 0 \end{pmatrix} \oplus  
A 
\begin{pmatrix} 0 \\  0 \\ 0 \\ 0 \\ * \\ 0 \end{pmatrix} \simeq 
A^{\oplus 3}, $$
we have $V_{\rho_A}^{\oplus 2} / {\rm Ker}(\partial_1) 
\simeq V_{\rho_A}$ 
and so the assumption (A2) is satisfied. 
Now, let 
$$P_\rho := 
\begin{pmatrix} 
\frac{\sqrt{(s^2 + 1)^2 - 4}}{2} & 1 \\ 
0 & -\frac{\sqrt{(s^2 + 1)^2 - 4}}{2} 
\end{pmatrix} 
\in {\rm sl}_2(A). 
$$
Since 
$$
H_1(I_{\mu}, {\rm Ad}(\rho_A))
= 
A \cdot 
P_{\rho} 
, 
$$ 
we can take 
$v_0$ as 
$P_{\rho}$
and so 
the assumption (A1) is 
satisfied. 

Next, let us 
see 
the presentation matrix $D$ in (2.6): 
$$
D = 
\begin{pmatrix} 
 & & & & & \\ 
 & {\Large {\rm Ad}(\rho_A) \left( \frac{\partial r}{\partial g_1} \right) } & & & {\Large {\rm Ad}(\rho_A) \left( \frac{\partial r}{\partial g_2} \right) } & \\ 
 & & & & & \\ 
1 & \frac{\sqrt{(s^2 + 1)^2 - 4}}{2} & 0 & 0 & 0 & 0 \\ 
0 & 0 & 0 & 0 & 0 & 0 \\ 
0 & 0 & 0 & 0 & 0 & 0 
\end{pmatrix}. $$ 

\noindent
Then 
the 
greatest common divisor 
$L_{\mu}(s) \in \mathbb{C}[[ s ]]$
of 
all 3-minors of $D$ is given by 
$$
L_{\mu}(x) \doteq 
\sqrt{ \{ (s^2 + 1)^2 - 1\} \{(s^2 + 1)^2 - 5 \} } \doteq 
s 
\in \mathbb{C}[[ s ]], 
$$ 
and so 
$V_{\rho_A}^{\oplus 2} / {\rm Im}(\delta) 
\simeq A^{\oplus 3} 
\oplus A / s A$. 
Hence by Theorem 2.7, we have 
$${\rm Sel}^{\rm h}_{\mu}({\rm Ad}(\rho_A))
\simeq A / s A 
\simeq \mathbb{C}. 
$$

\noindent 
{\bf Example 2.9.}
Let $K$ be the figure-eight knot. 
Next, we consider the case of $\gamma$ being the preferred longitude $\lambda$ corresponding to 
$\mu = g_1$. 
Let 
$\rho_A \: G_{K} \to {\rm SL}_{2} ( A )$ 
be the representation 
given by the following: 
$$\rho_A(g_1) 
= \begin{pmatrix} 
\frac{x + \sqrt{x^{2} - 4}}{2} & 
1 \\ 
0 & 
\frac{x - \sqrt{x^{2} - 4}}{2} 
\end{pmatrix}, $$ 
$$\rho_A(g_2) 
= \begin{pmatrix} 
\frac{x + \sqrt{x^{2} - 4}}{2} & 
0 \\ 
- (x^{2} - y(x) - 2) & 
\frac{x - \sqrt{x^{2} - 4}}{2} 
\end{pmatrix}, $$ 
where 
$A = \C \left[ \left[ x - \sqrt{\frac{5}{2}} \right] \right] $, and 
$y(x) =  \frac{x^2 + 1 + \sqrt{(x^2 - 1) (x^2 - 5)}}{2}$. 
We can rewrite these matrices as 
$$\rho_A(g_1) 
= \begin{pmatrix} 
\frac{ 
\left( s + \sqrt{\frac{5}{2}} \right) + 
\sqrt{ \left( s + \sqrt{\frac{5}{2}} \right)^{2} - 4}
}{2} & 
1 \\ 
0 & 
\frac{ 
\left( s + \sqrt{\frac{5}{2}} \right) - 
\sqrt{ \left( s + \sqrt{\frac{5}{2}} \right)^{2} - 4}
}{2} 
\end{pmatrix}, $$ 
$$\rho_A(g_2) 
= \begin{pmatrix} 
\frac{ 
\left( s + \sqrt{\frac{5}{2}} \right) + 
\sqrt{ \left( s + \sqrt{\frac{5}{2}} \right)^{2} - 4}
}{2} & 
0 \\ 
- \left\{ \left( s + \sqrt{\frac{5}{2}} \right)^{2} - 
y\left( s + \sqrt{\frac{5}{2}} \right) - 2 \right\} & 
\frac{ 
\left( s + \sqrt{\frac{5}{2}} \right) - 
\sqrt{ \left( s + \sqrt{\frac{5}{2}} \right)^{2} - 4}
}{2}\end{pmatrix}, $$ 
where we use 
$s := x - \sqrt{\frac{5}{2}}$. 
Similarly as Example 2.8, 
since 
${\rm Ker}(\partial_1) \simeq 
A^{\oplus 3}, $
we have $V_{\rho_A}^{\oplus 2} / {\rm Ker}(\partial_1) 
\simeq V_{\rho_A}$ 
and so the assumption (A2) is satisfied. 

Now, 
let 
$$P_\rho := 
\begin{pmatrix} 
\frac{\sqrt{\left( s + \sqrt{\frac{5}{2}} \right)^2 - 4}}{2} & 1 \\ 0 & -\frac{\sqrt{\left( s + \sqrt{\frac{5}{2}} \right)^2 - 4}}{2} 
\end{pmatrix} 
\in {\rm sl}_2(A), 
$$
and let $\lambda$ be the longitude with respect to 
$\mu = g_1$, 
and 
$h_{\rho}^{(1)}(\widetilde{\lambda})$ the reference generator 
of $H^1(G_K, {\rm Ad}(\rho_A))$ 
with respect to 
a lift $\widetilde{\lambda}$ 
in the universal covering 
of $\lambda$. 
Then
we have 
$$
T_{\lambda} \cdot h_{\rho}^{(1)}(\widetilde{\lambda}) (\widetilde{\mu}) = 
T_{\mu} \cdot h_{\rho}^{(1)}(\widetilde{\mu}) (\widetilde{\mu}), 
$$ 
(cf. \cite[Proposition 3.18]{porti1997torsion} for 
cohomology
), 
where $\displaystyle T_{\mu} := \frac{1}{2} \sqrt{ \left\{ \left( s + \sqrt{\frac{5}{2}} \right)^2 - 1 \right\} \left\{ \left( s + \sqrt{\frac{5}{2}} \right)^2 - 5 \right\} }$, 
$T_{\lambda} := 5 - 2 \left( s + \sqrt{\frac{5}{2}} \right)^2 \in A$, 
and 
by \cite[p.69]{porti1997torsion} 
$$
h_{\rho}^{(1)}(\widetilde{\mu}) (\widetilde{\mu}) = 
\begin{pmatrix} 
\frac{\partial M(s)}{\partial s} \cdot \frac{1}{M(s)} & 0 \\ 
0 & -\frac{\partial M(s)}{\partial s} \cdot \frac{1}{M(s)} 
\end{pmatrix} = 
\begin{pmatrix} 
\frac{1}{\sqrt{\left( s + \sqrt{\frac{5}{2}} \right)^2 - 4}} & 0 \\ 0 & -\frac{1}{\sqrt{\left( s + \sqrt{\frac{5}{2}} \right)^2 - 4}} 
\end{pmatrix} 
\in {\rm sl}_2(A), 
$$
where $M(s) = \frac{\left( s + \sqrt{\frac{5}{2}} \right) + \sqrt{\left( s + \sqrt{\frac{5}{2}} \right)^{2} - 4}}{2}\in A$ is 
the positive eigenvalue of $\rho_A(\mu)$. 
So we have 
$$
\langle h_{\rho}^{(1)}(\widetilde{\lambda}), P_{\rho} \otimes \widetilde{\mu} \rangle 
= \frac{T_{\mu}}{T_{\lambda}} \cdot 
\langle h_{\rho}^{(1)}(\widetilde{\mu}), P_{\rho} \otimes \widetilde{\mu} \rangle
= \frac{T_{\mu}}{T_{\lambda}}. 
$$ 
Hence, 
it is natural 
to take $v_0$ as 
$\displaystyle \frac{T_{\lambda}}{T_{\mu}} \cdot P_{\rho}. 
$

Next, let us 
see 
the presentation matrix 
$$
\begin{pmatrix} 
 & & & & & \\ 
 & {\Large {\rm Ad}(\rho_A) \left( \frac{\partial r}{\partial g_1} \right) } & & & {\Large {\rm Ad}(\rho_A) \left( \frac{\partial r}{\partial g_2} \right) } & \\ 
 & & & & & \\ 
\frac{T_{\lambda}}{T_{\mu}} & 
\frac{T_{\lambda}}{T_{\mu}} \cdot 
\frac{\sqrt{ \left( s + \sqrt{\frac{5}{2}} \right)^2 - 4}}{2} & 
0 & 0 & 0 & 0 \\ 
0 & 0 & 0 & 0 & 0 & 0 \\ 
0 & 0 & 0 & 0 & 0 & 0 
\end{pmatrix}. $$ 

\noindent
Then 
the 
greatest common divisor 
$L_{\lambda}(x) \in \C \left[ \left[ s \right] \right]$
of 
all 3-minors of 
this matrix 
is given by 
$$
L_{\lambda}(x) \doteq 
T_{\lambda} 
\doteq s 
\in \C \left[ \left[ s \right] \right], 
$$ 
and so 
$V_{\rho_A}^{\oplus 2} / {\rm Im}(\delta) 
\simeq A^{\oplus 3} 
\oplus A / s A$. 
Hence by Theorem 2.7, we have 
$${\rm Sel}^{\rm h}_{\lambda}({\rm Ad}(\rho_A))
\simeq A / s A 
\simeq \mathbb{C}
$$

\noindent 
{\bf Example 2.10.}
Let $K$ be the knot $5_2$, 
whose knot group is given by 
$$
G_{K} = \langle 
g_1, g_2 \ | \ 
g_1 g_2 g_1^{-1} g_2^{-1} g_1 g_2 g_1 = g_2 g_1 g_2 g_1^{-1} g_2^{-1} g_1 g_2
\rangle. 
$$
Let 
$\rho_A \: G_{K} \to {\rm SL}_{2} ( A )$ 
be the representation 
given by the following: 
$$\rho_A(g_1) 
= \begin{pmatrix} 
\frac{x + \sqrt{x^{2} - 4}}{2} & 
1 \\ 
0 & 
\frac{x - \sqrt{x^{2} - 4}}{2} 
\end{pmatrix}, $$ 
$$\rho_A(g_2) 
= \begin{pmatrix} 
\frac{x + \sqrt{x^{2} - 4}}{2} & 
0 \\ 
- (x^{2} - y(x) - 2) & 
\frac{x - \sqrt{x^{2} - 4}}{2} 
\end{pmatrix}, $$

\noindent
where 
$A = \C[[ \sqrt{x-\beta} ]]$, 
$\beta = 2.546 \dots \in \C$ is the simple root 
satisfying the equation
$$
\beta^4 + 2 \beta^3 - 5 \beta^2  - 14 \beta - 7 = 0, 
$$ 
and 
$y = y(x) \in A$ is the unique solution 
satisfying the equation
$$
y^3 - (x^2 + 1)y^2 +(3x^2 -2)y - 2x^2 + 1 = 0 
$$ 
and
$$
y(\beta) = \xi, 
$$
where $\xi = 3.132 \dots \in \C$ is the simple root 
satisfying the equation 
$$
\xi^4 - 6 \xi^3 + 11 \xi^2  - 6 \xi - 1 = 0. 
$$ 
We can rewrite these matrices as 
$$\rho_A(g_1) 
= \begin{pmatrix} 
\frac{(s^2 + \beta) + \sqrt{(s^2 + \beta)^{2} - 4}}{2} & 
1 \\ 
0 & 
\frac{(s^2 + \beta) - \sqrt{(s^2 + \beta)^{2} - 4}}{2} 
\end{pmatrix}, $$ 
$$\rho_A(g_2) 
= \begin{pmatrix} 
\frac{(s^2 + \beta) + \sqrt{(s^2 + \beta)^{2} - 4}}{2} & 
0 \\ 
- \{ (s^2 + \beta)^{2} - y(s^2 + \beta) - 2 \} & 
\frac{(s^2 + \beta) - \sqrt{(s^2 + \beta)^{2} - 4}}{2} 
\end{pmatrix}, $$ 
where we use 
$s := \sqrt{x - \beta}$. 
Similarly as Example 2.8, 
since 
${\rm Ker}(\partial_1) \simeq A^{\oplus 3}$, 
we have $V_{\rho_A}^{\oplus 2} / {\rm Ker}(\partial_1) 
\simeq V_{\rho_A}$ 
and so the assumption (A2) is satisfied. 
Now, let 
$$P_\rho := 
\begin{pmatrix} 
\frac{\sqrt{(s^2 + \beta)^2 - 4}}{2} & 1 \\ 
0 & -\frac{\sqrt{(s^2 + \beta)^2 - 4}}{2} 
\end{pmatrix} 
\in {\rm sl}_2(A), 
$$
Since 
$$
H_1(I_{\mu}, {\rm Ad}(\rho_A))
= 
A \cdot 
P_{\rho} 
, 
$$ 
we can take 
$v_0$ as 
$P_{\rho}$
and so 
the assumption (A1) is 
satisfied. 
Then 
the 
greatest common divisor 
$L_{\mu}(s) \in \mathbb{C}[[ s ]]$
of 
all 3-minors of $D$ 
associated with the meridian 
is given by 
$$
L_{\mu}(x) \doteq 
\sqrt{ (s^2 + \beta)^4 - \beta^4 } \doteq s 
\in \mathbb{C}[[ s ]]. $$ 
and by similar argument as Example 2.9, 
the 
greatest common divisor 
$L_{\lambda}(s) \in \mathbb{C}[[ s ]]$
of 
all 3-minors of $D$ 
associated with the longitude 
is given by 
$$
L_{\lambda}(s) \doteq 
5 (s^2 + \beta)^4 y-10 (s^2 + \beta)^4-
5 (s^2 + \beta)^2 y^2-7 (s^2 + \beta)^2 y+
31 (s^2 + \beta)^2+ 7y^2-7 y-21 
\doteq 1 
\in \mathbb{C}[[ s ]]. 
$$
Hence by Theorem 2.7, we have 
$${\rm Sel}^{\rm h}_{\mu}({\rm Ad}(\rho_A))
\simeq A / s A 
\simeq \mathbb{C}, \ 
{\rm Sel}^{\rm h}_{\lambda}({\rm Ad}(\rho_A))
= 0. 
$$

\ 

\begin{center}
{\bf 3. 
Adjoint homological Selmer modules for holonomy
representations} \\
\end{center}

Let $K$ be a hyperbolic knot in $S^3$ and let 
$hol \: G_K \rightarrow {\rm PSL}_2(\mathbb{C})$ be the holonomy representation attached to 
the complete hyperbolic structure on $\Int(E_K)$. 
It is known that $hol$ can be lifted to a representation 
$\rho_{0} \: G_K \to \SL_{2}(\O_{F_0, S_0})$, 
where $\O_{F_0, S_0}$ 
denotes the ring of 
$S_0$-integers of a finite algebraic
number field $F_0$ and a finite set $S_0$ of finite primes of $F$, and that 
$\rho_0$ is irreducible over $\C$ 
(\cite[3]{culler1983varieties}, \cite[3.2]{maclachlan2003arithmetic}).

Let $X(G_K)$ be the character variety of $\SL_2(\C)$-representations of $G_K$ and 
let $X_0(G_K)$ be the irreducible component of $X(G_K)$ containing $\chi_0 := \tr(\rho_0)$. 
Since $X_0(G_K)$ is defined over a finite algebraic number field 
(\cite{long2010fields}), 
$\overline{\Q}$-rational points 
are dense in $X(G_K)$, 
where 
$\overline{\Q}$ is the field of all algebraic numbers. 
For example, Dehn surgery points, 
which yield closed hyperbolic $3$-manifolds, 
are such algebraic points. 
Let $\chi_{\rho}$ be a $\overline{\Q}$-rational point in $X_0(G_K)$ which
corresponds to, up to conjugation, 
a representation 
$\rho_{\O_{F,S}} \colon G_K \to \SL_2(\O_{F,S})$, 
$\chi_{ \rho_{\O_{F,S}} } = \tr(\rho_{\O_{F,S}})$, 
where $\O_{F,S}$ is the ring of $S$-integers of a finite algebraic number field $F$ and 
$S$ is a finite set of primes of $F$. 
We write $\rho_{\C}$ for the scalar extension of 
$\rho_{\O_{F,S}}$ 
to $\C$. 
We assume that $\rho_{\C}$ is irreducible.

The following theorem is shown by Porti: \\ 

\noindent
{\bf Theorem 3.1} 
(\cite[3.3]{porti1997torsion}). {\em 
Let the notations and assumptions be as above.
The following assertions hold.

\noindent 
{\rm (1)} 
Suppose that 
$
\rho_{\O_{F,S}} 
|_{D_K}$ 
is hyperbolic. 
Then $\rho_{\C}$ is $\gamma$-regular, namely, 
$\Sel^{h}_{\gamma}( \Ad( \rho_{\C} ) )$ 
$= \{ 0 \}$ 
for some $\gamma \in D_K$ if and only if 
$\chi_{
\rho_{\O_{F,S}} 
}$ 
is a smooth and reduced point of $X_0(G_K)$. 

\noindent 
{\rm (2)} 
Suppose that 
$
\rho_{\O_{F,S}} 
|_{D_K}$ 
is parabolic. 
Let $\gamma \in D_K$ with 
$\rho_{\O_{F,S}} 
(\gamma) \neq \pm I$. 
Then $\rho_{\C}$ is $\gamma$-regular if and only if 
$\chi_{
\rho_{\O_{F,S}} 
}$ 
is a reduced point of $X_0(G_K)$. 
In particular, 
$(\rho_{0})_{\C}$ is $\gamma$-regular for any $\gamma \neq 1$. 

In these cases, the evaluation map
$$X_0(G_K) \to \C; \ \chi \mapsto \chi(\gamma)$$
is an analytic isomorphism on a neighborhood of 
$\chi_{
\rho_{\O_{F,S}} 
}$. 
}

\ 

\noindent
{\bf Theorem 3.2}. {\em 
Let $\rho_{\O_{F,S}} \colon G_K \to \SL_2(\O_{F,S})$ be a representation, and 
assume that the scalar extension 
$\rho_{\C} \colon G_K \to \SL_2(\C)$ of $\rho_{\O_{F,S}}$ is 
a $\gamma$-regular 
irreducible representation. 
Then $\Sel^{h}_{\gamma}( \Ad( \rho_{\O_{F.S}} ) )$
is a finitely 
generated torsion $\O_{F,S}$-module. 
In particular, $\Sel^{h}_{\gamma}( \Ad( \rho_{0} ) )$ 
is a finitely generated torsion $\O_{F_0. S_0}$-module for any 
$\gamma \neq 1$. 
} \\

\noindent
{\em Proof.}  
Let $\rho_F : G_K \rightarrow {\rm SL}_2(F)$ be the representation obtained from 
$\rho_{{\cal O}_{F,S}}$ by the extension of scalar to $F$. 
Since ${\rm Sel}^{\rm h}_{\gamma}({\rm Ad}(\rho_{F}))$ is an $F$-vector space and 
${\rm Sel}^{\rm h}_{\gamma}({\rm Ad}(\rho_{F})) \otimes_F \mathbb{C} = 
{\rm Sel}^{\rm h}_{\gamma}({\rm Ad}(\rho_{\C})) = 0$, 
${\rm Sel}^{\rm h}_{\gamma}({\rm Ad}(\rho_{F})) = 0$. 
By Lemma 1.1, ${\rm Sel}^{\rm h}_{\gamma}({\rm Ad}(\rho_{{\cal O}_{F,S}}))$ is 
a finitely generated 
torsion ${\cal O}_{F,S}$-module.
\hfill $\Box$\\

\noindent 
{\bf Example 3.3.}
Let $K$ be the figure-eight knot, 
whose knot group is given by 
$$G_{K} = \langle g_1, g_2 \ | \ g_1 g_2^{-1} g_1^{-1} g_2 g_1 = g_2 g_1 g_2^{-1} g_1^{-1} g_2 \rangle. $$

It is known that 
a lifting of the holonomy representation attached to 
the complete hyperbolic structure on $S^3 \setminus K$ 
is given, up to conjugation, by the following (\cite{riley1984nonabelian}): 
$$\rho_{0}: G_{K} \to 
{\rm SL}_{2} (\O_{F_{0}}) ; \ \ 
\rho_0(g_1) 
= \begin{pmatrix} 1 & 1 \\ 0 & 1 \end{pmatrix}, \ 
\rho_0(g_2) 
= \begin{pmatrix} 1 & 0 \\ \frac{1 + \sqrt{-3}}{2} & 1 \end{pmatrix}, $$ 
where 
$\O_{F_{0}} = \Z[\frac{1 + \sqrt{-3}}{2}]$ 
is the ring of integers of 
$F_0 = \mathbb{Q}(\sqrt{-3})$. 

Since $g_1 \neq 1$ in $D_K$, 
by Theorem 3.2, 
${\rm Sel}^{\rm h}_{}({\rm Ad}(\rho_{0}))$ is a finitely generated 
torsion $\O_{F_{0}}$-module. 
It is easy to see 
(by putting $x = 2$ in Example 2.8)
that 
assumptions (A1) and (A2) are satisfied, 
and 
by the straightforward computation 
and 
Theorem 2.7, 
we have 
$$
{\rm Sel}^{\rm h}_{}({\rm Ad}(\rho_{0)}) )
\simeq {\cal O}_{F_{0}} / \sqrt{-3} \ {\cal O}_{F_{0}}. 
$$

\noindent
{\bf Remark 3.4. }
In view of the analogy with number theory (Iwasawa theory),
our homological Selmer modules are expected to have some relations with
Reidemeister torsions associated with adjoint representations. 
In \cite[Chapter 4]{porti1997torsion}, 
Porti introduced the torsion function $\T_{(E_{K}, \mu)}$ on $X_0(G_K)$. 
According to his concrete computation worked out for the figure eight knot, 
we have 
$\T_{(E_{K}, \mu)} (\rho_{0}) \ 
= \left(\pm 1 / 2 \sqrt{(x^2 - 1) (x^2 - 5)} \right)|_{x=2}
= \pm \sqrt{-3} / 2$. 
See also \cite[5.2.1]{porti2015reidemeister}. 
Similar result holds for knot $5_2$. 
It would be interesting to pursue
connections between our homological Selmer modules and Porti’s torsions.

\begin{center}
{\bf 4. 
Adjoint homological Selmer modules for 
universal deformations} \\
\end{center}

We keep the same notations 
and assumptions 
as before. 
In this Section, 
we focus on the case of 
universal deformations, 
which is introduced in 
\cite{morishita2017universal} and \cite{kitayama2018certain}. 
For this reason, 
we first consider 
adjoint homological Selmer modules for
residual representations. 

\ 

\noindent
{\bf Proposition 4.1.} {\it  
There are infinitely many maximal ideal $\mathfrak{p}$ of ${\cal O}_{F,S}$  such that $\overline{\rho} := \rho_{{\cal O}_{F,S}} \; \mbox{\rm mod}\; \mathfrak{p} : G_K \rightarrow {\rm SL}_2(\mathbb{F}_{\mathfrak{p}})$ is an absolutely irreducible representation, where $\mathbb{F}_{\mathfrak{p}} := {\cal O}_F/\mathfrak{p}$, and ${\rm Sel}^{\rm h}_{\gamma}({\rm Ad}(\overline{\rho})) = 0$.}\\
\\
{\it Proof}. Let $\p$ be any maximal ideal of ${\cal O}_{F,S}$. Assume that $\overline{\rho} := \rho_{{\cal O}_{F,S}} \; \mbox{\rm mod}\; \p$ is not absolutely irreducible. Then there is a finite extension $F'$ of $F$ and a maximal ideal $\p$ of the ring of integers ${\cal O}_{F'}$ of $F'$ lying over $\p$ such that the representation $\overline{\rho}'$ obtained from $\overline{\rho}$ by the scalar extension to $\mathbb{F}_{\p} := {\cal O}_{F'}/\p$ is not an irreducible representation of $G_K$.  It suffices to show that the representation $\rho_{F'}$ obtained from $\rho_F$ by the scalar extension to $F'$ is not irreducible, for it contradicts to the irreducibility of $\rho$. Since $\overline{\rho}'$ is assumed to be reducible, there is a non-trivial  $G_K$-invariant subspace $V$ of $\mathbb{F}_{\p}^{\oplus 2}$, namely, $\dim_{\mathbb{F}_{\p}} V = 1$. Let $\Lambda$ be the inverse image of $V$ under the natural homomorphism ${\cal O}_{F'}^{\oplus 2} \rightarrow \mathbb{F}_{\p}^{\oplus 2}$ and let $W := \Lambda \otimes_{{\cal O}_{F'}} F'$. It is easy to see that $W$ is a subrepresentation of $\rho_{F'}$ over $F'$. Let
$\Lambda_{\p} := \Lambda \otimes_{{\cal O}_{F'}} {\cal O}_{\p}$ and $V_{\p} := \Lambda_{\p} \otimes_{{\cal O}_{\p}} F'_{\p}$, where ${\cal O}_{\p}$  and  $F'_{\p}$ are $\p$-adic completions of ${\cal O}_{F'}$ and $F'$ at $\p$, respectively. Then $\Lambda_{\p}$ is a free ${\cal O}_{\p}$-submodule of ${\cal O}_{\p}^{\oplus 2}$ and we have
$$ \dim_{F'} W = \dim_{F'_{\p}} V_{\p} = {\rm rank}_{{\cal O}_{\p}} \Lambda_{\p} = \dim_{\mathbb{F}_{\p}} V  =1.$$
Therefore $W$ is a non-trivial subrepresentation of $\rho_{F'}$. Hence we showed that  $\overline{\rho} := \rho_{{\cal O}_{F,S}} \; \mbox{\rm mod}\; \p$ is absolutely irreducible for any maximal ideal $\p$ of ${\cal O}_{F,S}$. 

On the other hand, by Theorem 3.2, 
there are only finitely many maximal ideal $\p$ of  ${\cal O}_{F,S}$ such that ${\rm Sel}^{\rm h}_{\gamma}({\rm Ad}(\rho_{{\cal O}_{F,S}})) \otimes_{{\cal O}_{F,S}} \mathbb{F}_{\p} \neq 0$.  Take a sufficiently large finite set $T$ of finite primes of $F$ containing  $S$ such that  ${\cal O}_{F,T}$ is a principal ideal domain and ${\rm Sel}^{\rm h}_{\gamma}({\rm Ad}(\rho_{{\cal O}_{F,T}})) \otimes_{{\cal O}_{F,T}} \mathbb{F}_{\p} = 0$ for any maximal ideal $\p$ outside $T$. We note that maximal ideals of ${\cal O}_{F,T}$ correspond bijectively to maximal ideals of ${\cal O}_{F,S}$ outside $T$. Let $\rho_{{\cal O}_{F,T}}$ be the representation obtained from $\rho_{{\cal O}_{F,S}}$ by the scalar extension to ${\cal O}_{F,T}$. Let $G$ denote $I_K$ or $G_K$. Let $\p$ be any maximal ideal of ${\cal O}_{F,T}$ and $\mathbb{F}_{\p} := {\cal O}_{F,T}/\p$. Consider the chain complex of ${\cal O}_{F,T}$-modules 
$$ C_2(G, {\rm Ad}(\rho_{{\cal O}_{F,T}})) \stackrel{\partial_2}{\rightarrow} C_1(G, {\rm Ad}(\rho_{{\cal O}_{F,T}})) \stackrel{\partial_1}{\rightarrow} C_0(G, {\rm Ad}(\rho_{{\cal O}_{F,T}})) 
\leqno{(4.1.1)}$$

\noindent
and the chain complex of $\mathbb{F}_{\p}$-modules
$$ C_2(G, {\rm Ad}(\overline{\rho})) \stackrel{\overline{\partial}_2}{\rightarrow} C_1(G, {\rm Ad}(\overline{\rho})) \stackrel{\overline{\partial}_1}{\rightarrow} C_0(G, {\rm Ad}(\overline{\rho})), \leqno{(4.1.2)}$$

\noindent
where $\overline{\rho} := \rho_{{\cal O}_{F,T}}$ mod $\p$ and $\overline{\partial}_i := \partial_i$ mod $\p$. 
From (4.1.1) we have  the exact sequence 
$$ 0 \rightarrow {\rm Im}(\partial_2) \rightarrow {\rm Ker}(\partial_1) \rightarrow H_1(G,{\rm Ad}(\rho_{{\cal O}_{F,T}})) \rightarrow 0.$$
By taking the tensor product over ${\cal O}_{F,T}$ with $\mathbb{F}_{\p}$,  we have
the exact sequence
$${\rm Im}(\partial_2)\otimes_{{\cal O}_{F,T}} \mathbb{F}_{\p} \rightarrow {\rm Ker}(\partial_1) \otimes_{{\cal O}_{F,T}} \mathbb{F}_{\p} \rightarrow H_1(G,{\rm Ad}(\rho_{{\cal O}_{F,T}})) \otimes_{{\cal O}_{F,T}} \mathbb{F}_{\p} \rightarrow 0.$$
Since ${\cal O}_{F,T}$ is a principal ideal domain, ${\rm Im}(\partial_i)$ and ${\rm Ker}(\partial_i)$ are free and so we have 
$${\rm Im}(\partial_2)\otimes_{{\cal O}_{F,T}} \mathbb{F}_{\p} = {\rm Im}(\overline{\partial}_2), \; {\rm Ker}(\partial_1) \otimes_{{\cal O}_{F,T}} \mathbb{F}_{\p} = {\rm Ker}(\overline{\partial}_1).$$
So, comparing with (4.1.2), we have 
$$H_1(G,{\rm Ad}(\rho_{{\cal O}_{F,T}})) \otimes_{{\cal O}_{F,T}} \mathbb{F}_{\p}= H_1(G, {\rm Ad}(\overline{\rho})).  \leqno{(4.1.3)}$$

\noindent
Consider the exact sequence of ${\cal O}_{F,T}$-modules 
 $$ H_1(I_K, {\rm Ad}(\rho_{{\cal O}_{F,T}})) \rightarrow 
 H_1(G_K,{\rm Ad}(\rho_{{\cal O}_{F,T}})) \rightarrow {\rm Sel}^{\rm h}_{\gamma}({\rm Ad}(\rho_{{\cal O}_{F,T}})) \rightarrow 0.$$
Taking the tensor product over ${\cal O}_{F,T}$ with $\mathbb{F}_{\p}$ and by (4.1.3), we have the exact sequence 
$$   H_1(I_K, {\rm Ad}(\overline{\rho})) \rightarrow H_1(G_K,{\rm Ad}(\overline{\rho})) \rightarrow {\rm Sel}^{\rm h}_{\gamma}({\rm Ad}(\rho_{{\cal O}_{F,T}})) \otimes_{{\cal O}_{F,T}} \mathbb{F}_{\p} \rightarrow 0$$
and hence
$$ {\rm Sel}^{\rm h}({\rm Ad}(\rho_{{\cal O}_{F,T}})) \otimes_{{\cal O}_{F,T}} \mathbb{F}_{\p} =  {\rm Sel}^{\rm h}_{\gamma}({\rm Ad}(\overline{\rho})).$$
By the choice of $T$, we have ${\rm Sel}^{\rm h}_{\gamma}({\rm Ad}(\overline{\rho})) = 0$.  \hfill $\Box$

\

Here, let us recall the universal deformation for knot group representations. 
For a local ring $A$, we denote by $\mathfrak{m}_A$ the maximal ideal of $A$. 
Let 
$k$ be a field with characteristic ${\rm char}(k) \neq 2$, 
and 
$\mathcal{O}$ a complete discrete valuation ring 
with $\mathcal{O}/\mathfrak{m_{\mathcal{O}}} = k$. 
Let $\overline{\rho}: G_K \rightarrow {\rm SL}_2(k)$ be a representation. 
The pair 
$(R, \rho)$ is called a {\em deformation} of $\overline{\rho}$ 
when $R$ is a complete local ${\cal O}$-algebra with $R/\mathfrak{m}_R = k$, 
and 
$\rho: G_K \rightarrow {\rm SL}_2(R)$ is a representation 
such that 
$\rho \ {\rm mod} \ \mathfrak{m}_R = \overline{\rho}$. 
Moreover, 
the pair $({\boldsymbol R}_{\overline{\rho}}, {\boldsymbol \rho})$ is called a {\em universal deformation} of $\overline{\rho}$ 
if 
$({\boldsymbol R}_{\overline{\rho}}, {\boldsymbol \rho})$ is a deformation of $\overline{\rho}$, 
and 
for all deformation $(R, \rho)$ of $\overline{\rho}$, 
there exists $\psi : {\boldsymbol R}_{\overline{\rho}} \rightarrow R$ 
such that  $\psi \circ {\boldsymbol \rho} \approx \rho$. 
Here, $\rho_1 \approx \rho_2$ 
means  
there exists $U \in I_2 + {\rm M}_2(\mathfrak{m}_R)$ 
such that $\rho_2(g) = U \rho_1(g) U^{-1}$ for all $g \in G_K$. 

\[\xymatrix{
\ 
& {\rm SL}_2({\boldsymbol R}_{\overline{\rho}}) \ar@{.>}[d]^{{}^{\exists 1} \psi} 
\ar@/{}^{5pc}/[dd]^{{\rm mod} \ \mathfrak{m}_{ {\boldsymbol R}_{\overline{\rho}} }} \\ 
G_K \ar[ur]^{\boldsymbol \rho} \ar[r]^{\rho} \ar[dr]_{\overline{\rho}} 
& {\rm SL}_2(R) \ar[d]^{{\rm mod} \ \mathfrak{m}_R} \\ 
\ 
& {\rm SL}_2(k) \\
}\]

\noindent
By the universal property, a universal deformation $({\boldsymbol R}_{\overline{\rho}}, {\boldsymbol \rho})$ 
of $\overline{\rho}$ is unique (if it exists) up to strict equivalence, 
and we call 
${\boldsymbol R}_{\overline{\rho}}$ 
the {\em universal deformation ring} of $\overline{\rho}$. 

Now, 
let $\p$ be a maximal ideal of ${\cal O}_F$ satisfying the properties in 
Proposition 4.1. 
Let ${\cal O}_{\p}$ be the $\p$-adic completion of ${\cal O}_F$.  
By \cite[Theorem 2.2.2]{morishita2017universal}, 
we have the universal deformation ${\boldsymbol \rho} \colon G_K \longrightarrow {\rm SL}_2({\boldsymbol R}_{\overline{\rho}})$ 
of $\overline{\rho}$. 
We assume that the universal deformation ring ${\boldsymbol R}_{\overline{\rho}}$ is a 
unique factorization
domain. 

\

\noindent
{\bf Theorem 4.2.} {\em For $\gamma \in D_K$ with $\gamma \neq 1$, the adjoint homological $\gamma$-Selmer module ${\rm Sel}^{\rm h}_{\gamma}({\rm Ad}({\boldsymbol \rho}))$ 
is a finitely generated 
torsion ${\boldsymbol R}_{\overline{\rho}}$-module.}\\
\\
{\em Proof.} 
By assumptions and direct calculations,  
we have 
$\dim H_0(E_K) = 0$ and $\dim H_0(\gamma) = 1$ 
over $\F_{\p}$ and 
the quotient field 
${Q( {\boldsymbol R}_{\overline{\rho}} ) }$ 
of ${\boldsymbol R}_{\overline{\rho}}$. 
By the exact sequence of relative homology, 
we have 
$$
0 \to {\rm Sel}_{\gamma}^h \to H_1(E_K, \gamma) \to H_0(\gamma) \to H_0(E_K) \to 0 
\leqno{(4.2.1)}
$$

\noindent
over $\F_{\p}$ and 
${Q( {\boldsymbol R}_{\overline{\rho}} ) }$. 
By Proposition 4.1, 
since ${\rm Sel}_{\gamma}^h (\Ad (\rhobar)) = 0$, 
we have 
$$
\dim H_1(E_K, \gamma, \Ad (\rhobar)) = 
\dim H_0(\gamma, \Ad (\rhobar)) - 
\dim H_0(E_K, \Ad (\rhobar)) = 
1, 
$$ 
and by the Euler characteristic argument, we have 
$$
\dim H_2(E_K, \gamma, \Ad (\rhobar)) = 
\dim {\it H}_1(E_K, \gamma, \Ad (\rhobar)) = 
1. 
\leqno{(4.2.2)}
$$ 

Next, 
let $(W, \gamma')$ be the 
2-dimensional CW complex 
$W$ (also defined in Section 2) 
with one vertex $\gamma'$, 
which is homotopically equivalent with $(E_K, \gamma)$. 
Consider 
the boundary map 
$\partial_2^{\prime} : C_2(W, \gamma' ) \to C_1(W, \gamma' )$
of the chain complex 
over $\F_{\p}$ and 
${Q( {\boldsymbol R}_{\overline{\rho}} ) }$. 
By (4.2.2), 
we have 
$${\rm corank}_{\F_{\p}} (\partial_2^{\prime}) = 
{\rm corank}_{Q( {\boldsymbol R}_{\overline{\rho}} ) } (\partial_2^{\prime}) = 0 {\rm \ or \ } 1. $$

Hence, we have 
$$\dim H_2(E_K, \gamma, \Ad (\rho_{Q( {\boldsymbol R}_{\overline{\rho}} ) }) ) = 0 {\rm \ or \ } 1, $$ 
and by the Euler characteristic argument, we have 
$$
\dim H_1(E_K, \gamma, \Ad (\rho_{Q( {\boldsymbol R}_{\overline{\rho}} ) }) ) = 
\dim H_2(E_K, \gamma, \Ad (\rho_{Q( {\boldsymbol R}_{\overline{\rho}} ) }) ) = 
0 {\rm \ or \ } 1. 
$$
By (4.2.1), 
we have 

$$
\dim H_1(E_K, \gamma, \Ad (\rho_{Q( {\boldsymbol R}_{\overline{\rho}} ) } ) = 
\dim H_0(\gamma, \Ad (\rho_{Q( {\boldsymbol R}_{\overline{\rho}} ) } ) - 
\dim H_0(E_K, \Ad (\rho_{Q( {\boldsymbol R}_{\overline{\rho}} ) } ) = 
1, 
$$ 
and so 
${\rm Sel}_{\gamma}^h (\Ad (\rho_{Q( {\boldsymbol R}_{\overline{\rho}} ) }) ) = 0$. 
Hence, 
${\rm Sel}^{\rm h}_{\gamma}({\rm Ad}({\boldsymbol \rho}))$ 
is a finitely generated 
torsion ${\boldsymbol R}_{\overline{\rho}}$-module. 
\hfill $\Box$

\ 

\noindent
{\bf Remark 4.3. }
The representation 
$\rho_A \: G_{K} \to {\rm SL}_{2} ( \C[[ \sqrt{x-1} ]] )$ 
in Example 2.8 
is a candidate of 
the universal representation 
by regarding the representation over the power series $\Z_p[[ \sqrt{x-1} ]]$ of 
the $p$-adic integers
for some prime number $p$. 
However, 
there still remains a rigorous discussion 
for the case of power series having square root in the variable. 

One way to avoid 
this problem 
is to replace the parameter $x = \tr(\rho_A(g_1))$ of the character variety 
by another parameter, 
such as $y = \tr(\rho_A(g_1 g_2))$, 
so that 
we may regard $\rho_A$ 
as the representation over the power series $\Z_p[[ y-1 ]]$, 
and we obtain 
$${\rm Sel}^{\rm h}_{\mu}({\rm Ad}(\rho_A))
\simeq \Z_p[[ y-1 ]] / (y-1) \Z_p[[ y-1 ]]. $$
See Example 5.3, where such kind of problem does not occur. 

\begin{center}
{\bf 5. 
Two-variable 
adjoint homological
Selmer modules for 
universal deformations} \\
\end{center}

Next, let us consider 
two-variable adjoint homological Selmer modules for universal deformations. 
Note that in our situation, 
``two-variable" means 
the parameters of 
the universal deformation and 
the character variety of 
one-dimensional representations 
of knot groups, 
namely the parameter $t$ of the Alexander polynomial. 
Let $A$ be an integral domain. 
Let $\rho \: G_K \rightarrow {\rm SL}_2(A)$ be a representation of $G_K$. 
We define the 
two-variable 
chain complex associated with $\rho$ 
for the knot complement $X_K$ 
as follows. 
Let $\widetilde{X} \to X_K$ be the universal cover of $X_K$. 
Let 
$\alpha : G_K \rightarrow G_K^{\rm ab} \simeq \Z 
= \langle t \rangle$ 
be the abelianization homomorphism. 
Then 
$V_{\rho}[t^{\pm 1}] = A[t^{\pm 1}] \otimes_A V_{\rho}$ 
is a right $A[G_K]$-module via 
$$
(a(t) \otimes {\boldsymbol v}).g 
:= a(t) \cdot \alpha(g) \otimes {\boldsymbol v} . \rho(g), 
$$ 
where $a(t) \in A[t^{\pm 1}]$, 
${\boldsymbol v} \in V_{\rho}$ and $g \in G_K$. 

We define the 
{\it two-variable $\rho$-twisted chain complex} 
$C_{*}(X_K; V_{\rho} [t^{\pm 1}] )$ 
{\it of} $X_K$ 
by 
$$
C_{*}(X_K; V_{\rho} [t^{\pm 1}] ) := 
V_{\rho}[t^{\pm 1}] \otimes_{A[G_K]}  C_{*}(\widetilde{X}), \ 
$$
and the 
{\it 
two-variable 
$\rho$-twisted homology 
}
$H_{i}(X_K; V_{\rho} [t^{\pm 1}] )$ by 
$$H_{i}(X_K; V_{\rho} [t^{\pm 1}] ) 
:= H_{i}(C_{*}(X_K; V_{\rho} [t^{\pm 1}] ) ). $$
Since 
the knot complement $X_K$ is the
Eilenberg-MacLane space $K(G_K, 1)$, 
we denote $H_{i}(X_K; V_{\rho} [t^{\pm 1}] )$ by 
$H_{i}(G_K; V_{\rho} [t^{\pm 1}] )$. 

For the case of the group $I_{\gamma}$, 
comparing with (2.2), 
the set 
$\{ v(t) \in V_{\rho} [t^{\pm 1}] \ | \ {\rho}(g_1) v(t) = t \cdot v(t) {\rho}(g_1)\}$
vanishes. 
Hence, 
we 
define the {\em two-variable adjoint homological $\gamma$-Selmer module}
attached to the
representation $\rho$ by the $A[t^{\pm 1}]$-module 
$${\rm Sel}^{\rm h}_{\gamma}({\rm Ad}(\rho)[t^{\pm 1}] ) := 
H_{1}(G_K; V_{\rho} [t^{\pm 1}] ). $$ 

\noindent
For the finitely generated 
torsion-ness of ${\rm Sel}^{\rm h}_{\gamma}({\rm Ad}(\rho)[t^{\pm 1}] )$, 
similar Euler characteritic arguments hold 
as \cite[Theorem 3.2.4]{kitayama2018certain}. 
Note that the condition 
$\det (t \cdot \boldrho(g) - I) \neq 0$ for some $g \in G_K$ and 
$\Delta_K(\boldrho; t) \neq 0$ 
always holds for the universal deformation $\boldrho$. 

\ 

\noindent
{\bf Theorem 5.1.} 
Let $({\boldsymbol R}_{\overline{\rho}}, \boldrho)$ be the universal deformation and 
assume that the universal deformation ring ${\boldsymbol R}_{\overline{\rho}}$ is a 
Noetherian 
unique factorization 
domain. 
For $\gamma \in D_K$ with $\gamma \neq 1$, 
the adjoint homological $\gamma$-Selmer module 
${\rm Sel}^{\rm h}_{\gamma}({\rm Ad}({\boldsymbol \rho}) [t^{\pm 1}] )$ 
is a finitely generated 
torsion ${\boldsymbol R}_{\overline{\rho}} [ t^{\pm 1} ]$-module. 

\ 

From now on, 
we keep the same notations as in Section 2. 
We study 
the $A[t^{\pm 1}]$-module structure
of the adjoint homological Selmer module
${\rm Sel}^{\rm h}({\rm Ad}(\rho)[t^{\pm 1}])$. 

We take 
again 
a Wirtinger presentation of $G_K$: 
$$G_K 
= \langle g_1,\dots , g_n \ | \ r_1 = \cdots = r_{n-1} = 1 \rangle, 
$$
where $g_1,\dots , g_n$ $(n \geq 2)$ represent 
meridians 
of a knot $K$, 
and let $\rho \: G_K \to \SL_2(A)$ be a representation. 
Let $V_{\rho}[t^{\pm 1}]$ be the representation space 
$$V_{\rho}[t^{\pm 1}] := 
{\rm sl}_2(A[t^{\pm 1}]) 
= A[t^{\pm 1}] v_1 
\oplus A[t^{\pm 1}] v_2 
\oplus A[t^{\pm 1}] v_3, $$ 
where 
$$v_1 := \begin{pmatrix} 0 & 1 \\ 0 & 0 \end{pmatrix}, \; 
v_2 := \begin{pmatrix} 1 & 0 \\ 0 & -1 \end{pmatrix}, \; 
v_3 := \begin{pmatrix} 0 & 0 \\ 1 & 0 \end{pmatrix}. $$
In the following, 
taking $\gamma = g_1$, 
we shall compute 
a presentation matrix 
of 
${\rm Sel}^{\rm h}_{}({\rm Ad}(\rho)[t^{\pm 1}]) = 
{\rm Sel}^{\rm h}_{\gamma}({\rm Ad}(\rho)[t^{\pm 1}])$ 
over $A[t^{\pm 1}]$, 
under some (mild) assumptions.  

Similar discussion as 
\cite[Section 3]{kitayama2018certain} 
holds  
by replacing 
the coefficient $V_{\rho}$ to $V_{\rho}[t^{\pm 1}]$, and 
the boundary maps 
$\partial_2$ 
to 
$
\partial_2[t^{\pm 1}] := 
\begin{pmatrix}
{\rm Ad}({\rho}) \circ \Psi \left( \frac{\partial r_i}{\partial g_j} \right)
\end{pmatrix}, 
$
where 
$
\Psi 
:= (\rho \otimes \alpha) \circ \pi : A[F] \to {\rm M}_2(A[t^{\pm 1}])$ 
is an $A$-algebra homomorphism, 
and $\alpha : G_K \to \Z$ is an abelianization. 

Similarly as 
Theorem 2.7, 
we have 
the following theorem. \\ 

\noindent 
{\bf Theorem 5.2.} {\em 
Let the notations and assumptions be as above. 
Then 
a presentation matrix of 
${\rm Sel}^{\rm h}({\rm Ad}({\rho}) [t^{\pm 1}] )$ 
is given by the matrix 
$\partial_2[t^{\pm 1}]$, and 
for an integer $d \geq 0$, 
the $d$-th Fitting ideal of ${\rm Sel}^h(\Ad(\rho) [t^{\pm 1}] )$ over $A[t^{\pm 1}]$ is 
generated by the greatest common divisor of all $(d + 3)$-minors of $\partial_2[t^{\pm 1}]$. 
} \\ 

We give herewith a concrete example, which provides a non-trivial result. 

\ 

\noindent 
{\bf Example 5.3.} 
Let $K$ be the figure-eight knot, 
whose knot group is given by 
$$G_{K} = \langle g_1, g_2 \ | \ g_1 g_2^{-1} g_1^{-1} g_2 g_1 = g_2 g_1 g_2^{-1} g_1^{-1} g_2 \rangle. $$
Let $\rho \: G_K \to \SL_2(\C)$ 
be the representation 
$$\rho(g_1) 
= \begin{pmatrix} 
\frac{\sqrt{\frac{5}{2}} + \sqrt{- \frac{3}{2}} }{2} & 
1 \\ 
0 & 
\frac{\sqrt{\frac{5}{2}} - \sqrt{- \frac{3}{2}} }{2}
\end{pmatrix}, $$ 
$$\rho(g_2) 
= \begin{pmatrix} 
\frac{\sqrt{\frac{5}{2}} + \sqrt{- \frac{3}{2}} }{2} & 
0 \\ 
\frac{5 + \sqrt{-15}}{4} & 
\frac{\sqrt{\frac{5}{2}} - \sqrt{- \frac{3}{2}} }{2}
\end{pmatrix}. $$
Consider a residual representation 
$\overline{\rho} \: G_K \to \SL_2(\F_{53})$ induced by 
$\rho$ 
such that 
$$\overline{\rho}(g_1) 
= \begin{pmatrix} 
19 & 
1 \\ 
0 & 
14
\end{pmatrix}, $$ 
$$\overline{\rho}(g_2) 
= \begin{pmatrix} 
19 & 
0 \\ 
44 & 
14
\end{pmatrix}. $$
Then 
we have the universal deformation 
$
\rho_R \: G_{K} \to {\rm SL}_{2} ( R )$ 
of $\overline{\rho}$ 
given by 
the following: 
$$\rho_R(g_1) 
= \begin{pmatrix} 
\frac{x + \sqrt{x^{2} - 4}}{2} & 
1 \\ 
0 & 
\frac{x - \sqrt{x^{2} - 4}}{2} 
\end{pmatrix}, $$ 
$$\rho_R(g_2) 
= \begin{pmatrix} 
\frac{x + \sqrt{x^{2} - 4}}{2} & 
0 \\ 
- ( x^2 - y(x) - 2 ) & 
\frac{x - \sqrt{x^{2} - 4}}{2} 
\end{pmatrix}, $$ 
where 
$R = \Z_{53} \left[ \left[ x - \sqrt{\frac{5}{2}} \right] \right]$, and 
$y( x ) = \frac{x^2 + 1 + \sqrt{(x^2 - 1) (x^2 - 5)}}{2}$. 
Similarly as Example 2.9, 
we can rewrite these matrices as 
$$\rho_A(g_1) 
= \begin{pmatrix} 
\frac{ 
\left( s + \sqrt{\frac{5}{2}} \right) + 
\sqrt{ \left( s + \sqrt{\frac{5}{2}} \right)^{2} - 4}
}{2} & 
1 \\ 
0 & 
\frac{ 
\left( s + \sqrt{\frac{5}{2}} \right) - 
\sqrt{ \left( s + \sqrt{\frac{5}{2}} \right)^{2} - 4}
}{2} 
\end{pmatrix}, $$ 
$$\rho_A(g_2) 
= \begin{pmatrix} 
\frac{ 
\left( s + \sqrt{\frac{5}{2}} \right) + 
\sqrt{ \left( s + \sqrt{\frac{5}{2}} \right)^{2} - 4}
}{2} & 
0 \\ 
- \left\{ \left( s + \sqrt{\frac{5}{2}} \right)^{2} - 
y\left( s + \sqrt{\frac{5}{2}} \right) - 2 \right\} & 
\frac{ 
\left( s + \sqrt{\frac{5}{2}} \right) - 
\sqrt{ \left( s + \sqrt{\frac{5}{2}} \right)^{2} - 4}
}{2}\end{pmatrix}, $$ 
where we use 
$s := x - \sqrt{\frac{5}{2}}$. 

Next, let us 
see 
the presentation matrix $\partial_2[t^{\pm 1}]$ in Theorem 5.2. 
Then 
the 
greatest common divisor 
$\Phi(s, t) \in \mathbb{Z}_{53} \left[ \left[ s \right] \right] [t^{\pm 1}]$
of 
all 3-minors of $\partial_2[t^{\pm 1}]$ is given by 
$$
\Phi(s, t) = 
(t - 1) \left[ t^2 - \left\{ 2 \left( s + \sqrt{\frac{5}{2}} \right)^2 - 3 \right\} t + 1 \right] 
\in \mathbb{Z}_{53} \left[ \left[ s \right] \right] [t^{\pm 1}]. $$ 

\noindent
Hence, we have 
$${\rm Sel}^{\rm h}_{\lambda}({\rm Ad}(\rho_R) [t^{\pm 1}] )
\simeq R[t^{\pm 1}] / 
(t - 1) \left[ t^2 - \left\{ 2 \left( s + \sqrt{\frac{5}{2}} \right)^2 - 3 \right\} t + 1 \right] 
R[t^{\pm 1}]. $$
\noindent
Moreover, 
by Example 2.9, 
we have 
$$\left. \frac{\Phi(s, t)}{t - 1} \right|_{t=1} = L_{\lambda}( s ). $$

\ 

\noindent 
{\bf Example 5.4.} 
Let $K$ be the knot $5_2$, 
whose knot group is given by 
$$
G_{K} = \langle 
g_1, g_2 \ | \ 
g_1 g_2 g_1^{-1} g_2^{-1} g_1 g_2 g_1 = g_2 g_1 g_2 g_1^{-1} g_2^{-1} g_1 g_2
\rangle. 
$$
Consider a residual representation 
$\overline{\rho} \: G_K \to \SL_2(\F_{17})$ 
such that 
$$\overline{\rho}(g_1) 
= \begin{pmatrix} 
1 & 
1 \\ 
0 & 
1
\end{pmatrix}, $$ 
$$\overline{\rho}(g_2) 
= \begin{pmatrix} 
1 & 
0 \\ 
2 & 
1
\end{pmatrix}. $$
Then 
we have the universal deformation 
$
\rho_R \: G_{K} \to {\rm SL}_{2} ( R )$ 
of $\overline{\rho}$ 
given by 
the following: 
$$\rho_R(g_1) 
= \begin{pmatrix} 
\frac{x + \sqrt{x^{2} - 4}}{2} & 
1 \\ 
0 & 
\frac{x - \sqrt{x^{2} - 4}}{2} 
\end{pmatrix}, $$ 
$$\rho_R(g_2) 
= \begin{pmatrix} 
\frac{x + \sqrt{x^{2} - 4}}{2} & 
0 \\ 
- ( x^2 - y(x) - 2 ) & 
\frac{x - \sqrt{x^{2} - 4}}{2} 
\end{pmatrix}, $$ 
where 
$R = \Z_{17} \left[ \left[ x - \beta \right] \right]$, 
$\beta = \sqrt{-1} \cdot 1.00098 \dots \in \Z_{17}$ is the simple root 
satisfying the equation 
$$
20 \beta^6 - 126 \beta^4 + 196 \beta^2 + 343 = 0, 
$$ 
and 
$y = y(x) \in R$ is, 
by Hensel’s lemma, 
the unique solution satisfying the equation 
$$
y^3 - (x^2 + 1)y^2 +(3x^2 -2)y - 2x^2 + 1 = 0 
$$ 
and
$$
y(\beta) = \xi, 
$$
where $\xi = -2.493 \dots \in \Z_{17}$ is the simple root 
satisfying the equation 
$$
2 \xi^3 - 2 \xi^2 - 11 \xi + 16 = 0. 
$$ 

\noindent
We can rewrite these matrices as 
$$\rho_A(g_1) 
= \begin{pmatrix} 
\frac{ 
\left( s + \beta \right) + 
\sqrt{ \left( s + \beta \right)^{2} - 4}
}{2} & 
1 \\ 
0 & 
\frac{ 
\left( s + \beta \right) - 
\sqrt{ \left( s + \beta \right)^{2} - 4}
}{2} 
\end{pmatrix}, $$ 
$$\rho_A(g_2) 
= \begin{pmatrix} 
\frac{ 
\left( s + \beta \right) + 
\sqrt{ \left( s + \beta \right)^{2} - 4}
}{2} & 
0 \\ 
- \left\{ \left( s + \beta \right)^{2} - 
y\left( s + \beta \right) - 2 \right\} & 
\frac{ 
\left( s + \beta \right) - 
\sqrt{ \left( s + \beta \right)^{2} - 4}
}{2}\end{pmatrix}, $$ 
where we use 
$s := x - \beta$. 

Next, let us 
see 
the presentation matrix $\partial_2[t^{\pm 1}]$ in Theorem 5.2. 
Then 
the 
greatest common divisor 
$\Phi(x, t) \in \mathbb{Z}_{17} [ \left[ s \right] ] [t^{\pm 1}]$
of 
all 3-minors of $\partial_2[t^{\pm 1}]$ is given by 
\begin{align*}
\Phi( s, t) =
(t - 1) 
& [ 
\left\{ 2 ( s + \beta )^4 y-4 ( s + \beta )^4-2 ( s + \beta )^2 y^2-2 ( s + \beta )^2 y+10 ( s + \beta )^2+2 y^2-2 y-6 \right\} t^2  \\
&+ \left\{ ( s + \beta )^4 y-2 ( s + \beta )^4-( s + \beta )^2 y^2-3 ( s + \beta )^2 y+11 ( s + \beta )^2+3 y^2-3 y-9\right) t  \\ 
&+ \left\{ 2 ( s + \beta )^4 y-4 ( s + \beta )^4-2 ( s + \beta )^2 y^2-2 ( s + \beta )^2 y+10 ( s + \beta )^2+2 y^2-2 y-6 \right\} ] \\ 
&\in \mathbb{Z}_{17} [ \left[ s \right] ] [t^{\pm 1}]. 
\end{align*}
\noindent
Hence, we have 
$${\rm Sel}^{\rm h}_{\lambda}({\rm Ad}(\rho_R) [t^{\pm 1}] )
\simeq R[t^{\pm 1}] / \Phi( s, t) R[t^{\pm 1}]. $$
\noindent
Moreover, 
by Example 2.10, 
we have 
$$\left. \frac{\Phi(s, t)}{t - 1} \right|_{t=1} = L_{\lambda}(s). $$

\ 

Based on these examples, we have the following conjecture:

\

\noindent
{\bf Conjecture 5.5. }
Let $({\boldsymbol R}_{\overline{\rho}}, \boldrho)$ be the universal deformation and 
assume that the universal deformation ring ${\boldsymbol R}_{\overline{\rho}}$ is a 
Noetherian 
unique factorization domain. 
Let $\Phi(s, t) \in {\boldsymbol R}_{\overline{\rho}} [ t^{\pm 1} ]$ be 
the greatest common divisor of 
all 3-minors of $\partial_2[t^{\pm 1}]$ and 
$L_{\lambda}(s) \in {\boldsymbol R}_{\overline{\rho}}$ be 
the greatest common divisor of 
all 3-minors of $D$ in (2.6) associated with the longitude. 
Then we have 
$$\left. \frac{\Phi(s, t)}{t - 1} \right|_{t=1} = L_{\lambda}(s). $$

\

This conjecture is considered as a knot-theoretic analogue of Theorem $1$ for two-variable adjoint Selmer groups in Hida-Tilouine-Urban \cite{hida1997adjoint} (see also \cite[Theorem 6.3]{hida2000adjoint}, \cite[Chapter 5]{hida2006hilbert}). 
It would be interesting to pursue this analogy in detail. 

\pagebreak

\printbibliography

@article{culler1983varieties,
  title={Varieties of group representations and splittings of 3-manifolds},
  author={Culler, Marc and Shalen, Peter B.},
  journal={Annals of Mathematics},
  volume={117},
  pages={109--146},
  year={1983}
}

@article{fox1953free,
  title={Free differential calculus. {I}. {D}erivation in the free group ring},
  author={Fox, Ralph H.},
  journal={Annals of Mathematics},
  volume={57},
  number={2},
  pages={547--560},
  year={1953}
}

@incollection{greenberg1989iwasawa,
  title={Iwasawa Theory for $ p $-adic Representations},
  author={Greenberg, Ralph},
  booktitle={Algebraic Number Theory—in Honor of K. Iwasawa},
  volume={17}, 
  pages={97--137},
  year={1989},
  publisher={Advanced Studies in Pure Mathematics}
}

@article{hida2000adjoint,
  title={Adjoint Selmer groups as {I}wasawa modules},
  author={Hida, Haruzo},
  journal={Israel Journal of Mathematics},
  volume={120},
  number={2},
  pages={361--427},
  year={2000},
  publisher={Springer}
}

@book{hida2006hilbert,
  title={Hilbert {M}odular {F}orms and {I}wasawa {T}heory},
  author={Hida, Haruzo},
  year={2006},
  publisher={Oxford University Press}
}

@book{hida2000modular,
  title={Modular {F}orms and {G}alois {C}ohomology},
  author={Hida, Haruzo},
  year={2000},
  publisher={Cambridge University Press}
}

@article{hida1997adjoint,
  title={Adjoint modular {G}alois representations and their {S}elmer groups},
  author={Hida, Haruzo and Tilouine, Jacques and Urban, Eric},
  journal={The Proceedings of the National Academy of Sciences},
  volume={94},
  number={21},
  pages={11121--11124},
  year={1997},
  publisher={National Acad Sciences}
}

@book{kawauchi2012survey,
  title={A Survey of Knot Theory},
  author={Kawauchi, Akio},
  year={2012},
  publisher={Birkh{\"a}user}
}

@article{kitayama2018certain,
  title={On certain {L}-functions for deformations of knot group representations},
  author={Kitayama, Takahiro and Morishita, Masanori and Tange, Ryoto and Terashima, Yuji},
  journal={Transactions of the American Mathematical Society},
  volume={370},
  number={5},
  pages={3171--3195},
  year={2018}
}

@article{long2010fields,
  title={Fields of definition of canonical curves},
  author={Long, Darren D. and Reid, Alan W.},
  journal={Interactions between hyperbolic geometry, quantum topology and number theory},
  volume={541},
  pages={247--257},
  year={2010}
}

@book{maclachlan2003arithmetic,
  title={The arithmetic of hyperbolic 3-manifolds},
  author={Maclachlan, Colin and Reid, Alan W.},
  volume={219},
  year={2003},
  publisher={Springer}
}

@book{morishita2012knots,
  title={Knots and {P}rimes: {A}n {I}ntroduction to {A}rithmetic {T}opology},
  author={Morishita, Masanori},
  year={2012},
  publisher={Universitext, Springer, London}
}

@article{morishita2017universal,
  title={On the universal deformations for {SL}$_2$-representations of knot groups},
  author={Morishita, Masanori and Takakura, Yu and Terashima, Yuji and Ueki, Jun},
  journal={Tohoku Mathematical Journal},
  volume={69},
  number={1},
  pages={67--84},
  year={2017},
  publisher={Tohoku University, Mathematical Institute}
}

@article{porti2015reidemeister,
  title={Reidemeister torsion, hyperbolic three-manifolds, and character varieties},
  author={Porti, Joan},
  journal={Handbook of group actions, Advanced Lectures in Mathematics},
  volume={IV},
  number={41},
  pages={447--507},
  year={2018},
  publisher={Int. Press, Somerville, MA}
}

@book{porti1997torsion,
  title={Torsion de {R}eidemeister pour les vari{\'e}t{\'e}s hyperboliques},
  author={Porti, Joan},
  volume={612},
  year={1997},
  publisher={American Mathematical Society}
}

@article{riley1984nonabelian,
  title={Nonabelian representations of 2-bridge knot groups},
  author={Riley, Robert},
  journal={The Quarterly Journal of Mathematics},
  volume={35},
  number={2},
  pages={191--208},
  year={1984},
  publisher={Oxford University Press}
}

@article{tange2021non,
  title={Non-acyclic {SL}$_2$-representations of {T}wist {K}nots, $-3$-{D}ehn {S}urgeries, and $L$-functions},
  author={Tange, Ryoto and Tran, Anh T. and Ueki, Jun},
  journal={International Mathematics Research Notices},
  year={2021},
}

@article{thurston1979geometry,
  title={The geometry and topology of three-manifolds},
  author={Thurston, William P.},
  journal={Lecture Note, Princeton},
  year={1979},
  publisher={Princeton University Princeton, NJ}
}

\noindent
Takahiro Kitayama kitayama@ms.u-tokyo.ac.jp \\
Graduate School of Mathematical Sciences, 
the University of Tokyo, 
3-8-1, Komaba, Meguro-ku, Tokyo, 153-8914, Japan

\ 

\noindent
Masanori Morishita morishita.masanori.259@m.kyushu-u.ac.jp \\ 
Faculty of Mathematics, 
Kyushu University, 
744, Motooka, Nishi-ku, Fukuoka, 819-0395, Japan 

\ 

\noindent
Ryoto Tange rtange.math@gmail.com \\
Department of Mathematics, School of Education, Waseda University, 1-104, Totsuka-cho,
Shinjuku-ku, Tokyo, 169-8050, Japan 

\ 

\noindent
Yuji Terashima yujiterashima@tohoku.ac.jp \\ 
Graduate School of Science, 
Tohoku University, 
6-3, Aoba, Aramaki-aza, Aoba-ku, Sendai, 980-8578, Japan 

\end{document}